# Describe Prime number gaps pattern by Logistic mapping

Wang Liang
E-mail: wangliang.f@gmail.com

**[Abstract]** Based on symbolic dynamics methods, we show the primes gap pattern could be described by the chaos orbit of Logistic mapping $X(k+1)=1-uX(k)^2$, $u=1.5437$. If so, there will be arbitrarily many twin primes. This connection may provide a new avenue to prove twin primes conjecture.
**[Keywords]** Prime number, Symbolic dynamics, Logistic Mapping

## 1 Introduction

In this paper, we find:

**Theme 1**: Logistic mapping $x_{n+1} = 1 - u x_n^2$, $x_n \in [-1,1]$, $u \approx 1.5437$ could describe primes gap pattern.

The series orbits $x_{n+1} = 1 - u x_n^2, u < 2$ are shown in the bifurcation figure [Fig 1].

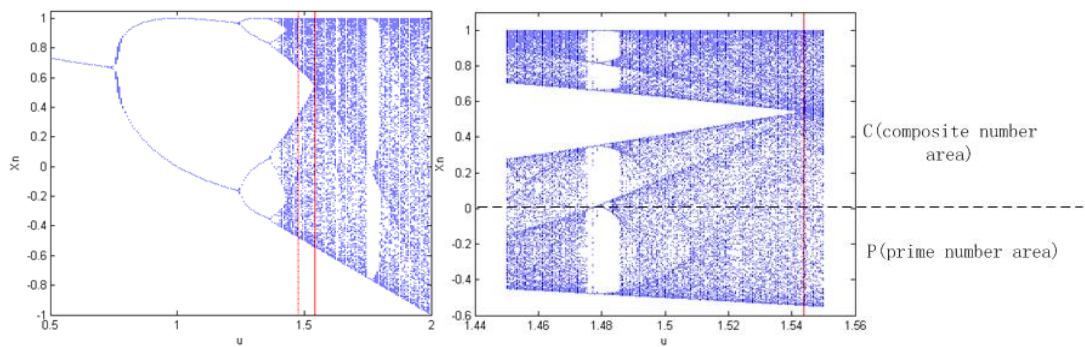

Fig1. Bifurcation figure of Logistic mapping, the red line in right figure corresponding to u=1.5437

$u = 1.5437$ in the bifurcation figure is the converging point of 'two bands' and 'single band'. The orbit of its corresponding Logistic mapping is chaos orbit.

Then we connect its orbit with prime gap based on symbolic dynamics methods. For orbit point x<0, we mark it "P", which represent the prime number point. For point x>0, mark it "C", represent composite number.

Because it's a chaos orbit, there will be arbitrarily many pattern "….PCP…." in its orbit, which means there will be arbitrarily many twin primes.

The statistics of difference of consecutive prime numbers are shown in Fig.2.A. Corresponding results of difference of "prime numbers" constructed by chaos orbit of Logistic mapping are shown in Fig.1.B. Their figure shapes are also very similar.

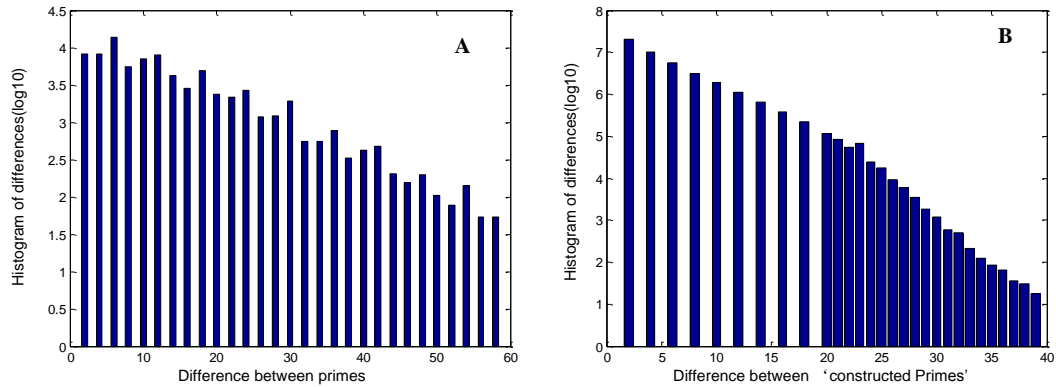

Fig2. Statistic of difference of consecutive primes (A) , and constructed primes(B) by logistic map

The following sections will introduce how to connect the prime gap with Logistic mapping.

## 2 Related work

The distribution of the prime numbers among the integers seems somewhat random, but it's surely produced by simple sieve method. This makes people believe the primes system is a chaos system. So some new methods in chaos and statistic theory are applied to study the pattern of prime numbers. For example in [1][2] the fractal character of primes was investigated, while in [3] the appropriately defined Lyapunov exponents for the distribution of primes were calculated numerically. Power-law behavior in the distribution of primes and correlations in prime numbers were found in [4]. The differences of consecutive prime numbers are also well studied. Besides the well known 3 period oscillations, the same special periodic behavior also appears in Dirichlet classification and increment (difference of difference) of consecutive prime [5][6]. Most of these researches all show interesting numerical feature in prime numbers, but unfortunately, they all can't deduce the strict theory results.

Obviously, if primes distribution is chaos, we need corresponding tools to study it. It's basic idea of this paper.

Now there have been many tools for analyzing chaos system, such as renormalization group, symbolic dynamics, etc. Here we select symbolic dynamics theory to analyze the prime gap pattern. Mainly because this theory also allows us to connect many complex systems with very simple mappings. Symbolic dynamic is originally a pure math theory, but its application in chaos analyzing makes it become a simple and powerful tool for dynamic system research, especially for discrete system [7].It has been successfully applied in some famous problems like 'Smale horseshoe mapping' and 'Three body problem'.

## 3 Brief introduction to symbolic dynamics theory

The symbolic dynamic theory is not very familiar for most researchers in number theory area, so we give a brief introduction for it first. Here we mainly refer to the "Applied symbolic dynamics theory":

Symbolic dynamics is a coarse-grained description of dynamics by taking into account the

"geometry". It could be easily explained by one-dimensional unimodal mappings $x_{n+1} = f(u, x_n)$ (Fig 2).

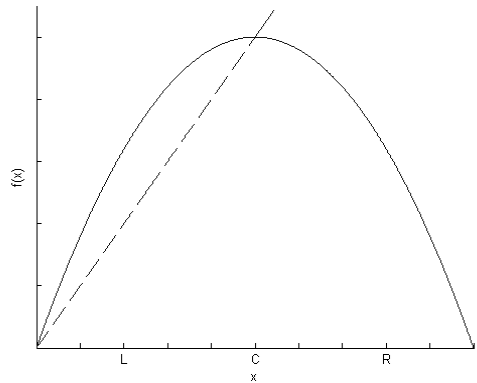

Fig 2. Symbol sequence in one dimension unimodal mapping

In Fig6, point 'c' in x axis corresponds to the maximal values of function $f(u, x)$. It divides the x axis into left part and right part. Through the iterated function $x_{n+1} = f(u, x_n)$ we could get an orbit from an initial point $x_0$: $x_0, x_1 = f(x_0), x_2 = f(x_1), \cdots, x_n = f(x_{n-1}) \cdots$

If we only care the relative place of $x_i$: the point on the left of 'c' is marked as 'L' and the right is 'R. So this orbit could be described by a symbolic sequence like $'LLLRRRC\cdots'$, here, $x_i = \begin{cases} L, x_i < c \\ C, x_i = c \\ R, x_i > c \end{cases}$.

Normally, we use basic sequence to represent the whole period sequence, for example $RL \leftrightarrow RLRLRL\cdots$

Our discussion mainly refers to three rules of symbolic dynamics:
(1) Kneading symbolic sequence. Any orbit can be described by a symbolic sequence, but some symbolic sequences couldn't correspond to an actual orbit. MSS table describe the basic 'admissible' sequence. Such sequence is also called kneading symbolic sequence, which describe the orbits beginning from the maximal value point of mapping. The kneading sequence must satisfy: $M < \delta(M)$, here $\delta$ is left shift function. How to compare two symbol sequence could be found in[7].
(2) The relation between symbolic sequence and parameter in a mapping. For a special mapping $x_{n+1} = f(u, x_n)$, we could get a parameter $u$ from one kneading symbolic sequence. That means we could always find the actual orbit corresponding to a kneading sequence.

Here we select Logistic mapping $x_{n+1} = 1 - ux_n^2$. The relation of its parameter u and related symbols is shown in Table 1, which could be found in some books for 'Applied symbolic dynamics' [9].

| u | 1.0 | 1.3107 | 1.3815 | 1.40115 | 1.4304 | 1.5437 | 1.754 | 2.0 |
|---|---|---|---|---|---|---|---|---|
| Symbols | $RC$ | $RLRC$ | $RLR^3LRC$ | $R^{*\infty}$ | $RLR^2(RL)^\infty$ | $RLR^\infty$ | $RLC$ | $RL^\infty$ |

Table 1. The relation of parameter u and related sequence

(3) For any two kneading sequence $M_1, M_2$, if $M_1 < M_2$, then for a special mapping, their corresponding parameter $u(M_1) < u(M_2)$.

## 4 symbolic dynamics model for primes

Because there seems still no related work about symbolic dynamics model for primes, we need build its model from basic sieve method.

Our research route is very clear. We use the symbol sequence to describe the sieve method. Then we could connect such sequences with a Logistic mapping. The dynamics characters of related mapping will represent pattern of primes. Here we have:

**Theme 2**: Primes gap pattern could be described by symbol kneading sequence $RLR^\infty$. Dynamics of this sequence could be describe by logistic mapping $x_{n+1} = 1 - ux_n^2$, $u \approx 1.5437$.

### 4.1 symbols expression for prime numbers

Eratosthenes Sieve method of two thousand years ago is still the main method to create the prime numbers. Given a list of the numbers between 1 and N, starting with 2, erase all multiples of 2 up to N, other than 2 itself. Call the remaining set $P_2$. Then return to the beginning and taking the first number greater than 2 and erase all of its multiples up to N, again other than the number itself. Repeating this operation, we could get $P_3, P_4, P_5, \cdots P_R$. If $R = N^{1/2}$, the set $P_R(N)$ will all be prime numbers. We could give a new representation for this method.

Here the nature number is regarded as the points in a line. The value of nature number is expressed by their position in this line. So 'erased' point is marked by one symbol and 'saving' point by other symbol. So we could use only two symbols to describe the sieve method. Here we select 'R,L'. In fact, any two symbols like '1,0' will also do. We first define the prime number by symbols and then describe the sieve method.

We select $LLLLLLLL\cdots = (L)^n$ to represent the all natural number points $(0,1,2,3,\cdots\cdots)$ in a line.

For prime number 2, all 'L' in the places that could be divided by 2 are 'erased' and replaced by 'R'. Here the first 'L' presenting 0 is also 'erased'. We get:

$RLRLRL\cdots = (RL)^n$.

For 3, all 'L' in the places that could be divided by 3 are replaced by 'R', We get:

$RLLRLLRLL\cdots = (RLL)^n$

Thorough this method, prime number could be expressed by a periodic oscillation symbol sequence. This operation is shown in Fig 3.

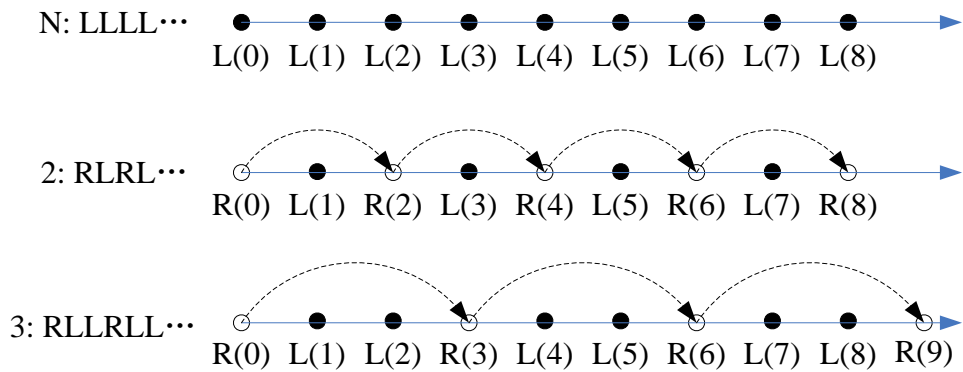

Fig 3. Period expression for prime number

By this method, for any prime number $p_i$, we could get its symbol expression $(RL^{p_i-1})^n$

In this paper, we use the basic symbolic sequence of one period to present the whole period sequence. So the prime number could be described as follows:

$2: M_2 = RL$;
$3: M_3 = RLL$;
$5: M_5 = RLLLL$;     $p_i$ is the ith prime number.
$\vdots$
$p_i: M_{P_i} = RL^{p_i-1}$.

### 4.2 Composition rule

After defining the symbolic expression for prime numbers, we could define their composition rule.

For example, we cut all the numbers that could be divided by 2 and 3. According to the Sieve rule, all 'L' in the place that could be divided by 2 or 3 is 'erased' and replaced by 'R'. We could get *RLRRRL*, it's a 6 period sequence. The composition of $M_2, M_3$ is very similar with dot product of vectors. This operation is shown in Fig 4.

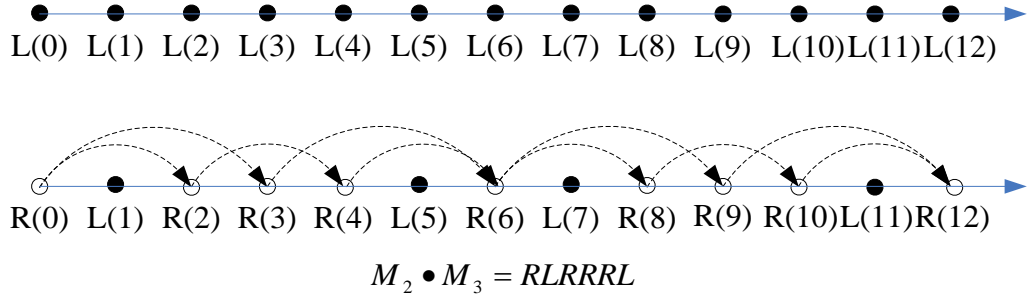

$$M_2 \bullet M_3 = RLRRRL$$

Fig 4 Composition of $M_2, M_3$

We use '$\bullet$' operator to describe composition operation. According to Fig4, denote the composition of $M_2, M_3$ is: $P_2 \bullet P_3 = RL \bullet RLL$. Here RL and RLL all represent the period sequence, so we could expand them as follows:

$$RL \bullet RLL \Leftrightarrow (RL)^3 \bullet (RLL)^2 = (RLRLRL) \bullet (RLLRLL)$$
$$= (R \bullet R)(L \bullet L)(R \bullet L)(L \bullet R)(R \bullet L)(L \bullet L) = RLRRRL$$

Here 'R' presents the 'erased' function. We could get the composition rules of two symbols:

1. $(R \bullet R)$ means 'erased' twice: $R \bullet R = R$.

2. $(L \bullet R)$ and $(R \bullet L)$ means 'erased' by 2 or 3: $L \bullet R = R \bullet L = R$.

3. $(L \bullet L)$ is the 'saving' point: $(L \bullet L) = L$

We could define composition rule for any two sequence::

**Definition 1**: '$\bullet$' composition rule is defined as follows:

If two symbol sequence have the equal length:

$$A = a_1 a_2 a_3 \cdots a_n, B = b_1 b_2 b_3 \cdots b_n$$
$$A \bullet B = (a_1 a_2 a_3 \cdots a_n) \bullet (b_1 b_2 b_3 \cdots b_n) = (a_1 \bullet b_1)(a_2 \bullet b_2)(a_3 \bullet b_3) \cdots (a_n \bullet b_n)$$
$$here, a, b \in \{R, L\}, a \bullet b = \begin{cases} L, a = b = L \\ R, a = b = R \\ R, a \neq b \end{cases}$$

If two symbols are different length, we just need extend their basic period sequence to same length:

$$A = a_1 a_2 a_3 \cdots a_n, B = b_1 b_2 b_3 \cdots b_m$$
$$A \bullet B = (a_1 a_2 a_3 \cdots a_n)^m \bullet (b_1 b_2 b_3 \cdots b_m)^n$$

## 4.3 Primes dynamics

We could use $D_i, D_1 = M_2, D_2 = M_2 \bullet M_3, \cdots, D_i = M_2 \bullet M_3 \bullet \cdots \bullet M_{pi}$ to describe the dynamical character of prime system. For example: $D_2 = M_2 \bullet M_3 = RLRRRL$ is 6 period. It contains the period oscillation 'information' of prime number 2 and 3. The dynamical character of whole prime system could be depicted as:

$$D_\infty = M_2 \bullet M_3 \bullet \cdots \bullet M_{pi} \bullet \cdots$$

In the dynamical expression of prime number, the length of consecutive R must be odd and there is no consecutive L. So to be a kneading sequence, it must satisfy:

**If the first length of consecutive R is n, the length of following consecutive R is no more than n.** There have been many researches for the difference of primes. Its proof is even more difficult than Riemann hypothesis.

Here we assume $D_2, D_3, \cdots, D_{pi}, \cdots$ are kneading sequence, which means all the consecutive difference of primes less than $p_1 \times p_2 \times p_3 \times \cdots \times p_i$ is not larger than $p_i$. Obviously, it's not true. We could easily find the counter example. But is part of $M_{pi}$ is kneading sequence? For example, only "L" point before $p_i^2$ are prime numbers. So we only need consider the first $p_i^2$ symbols. Here we assume:

**Theme 3: The first $p_i^2 + 1$ symbols of $D_i$ make up of a kneading sequence.**

This conjecture means all the consecutive difference of primes less than $p_i^2$ is not larger than $p_i$. Here we mark the g(N) as the consecutive difference of primes before N. We need:

$$g(N) < N^{1/2}$$

Now this is still a conjecture (Legendre's conjecture). Till now, the best result about upper bound of $g(N)$ is $g(N) < cN^{21/40}$ [8]. The precise upper bounds of prime gaps has no much impact for our discussion. We only need select a shorter sequence to ensure it's the kneading sequence.

**Theorem 4:** $D_1 < D_2 < D_3 < \cdots < D_i < D_{i+1} < \cdots < D_\infty$.

**Proof:**

For any $D_i$ and $D_{i+1}$, they could be written as:

$$D_i = \underbrace{RLRRR\cdots RL}_{A}\cdots \qquad D_{i+1} = \underbrace{RLRRR\cdots RR}_{B}\cdots$$

Here the length of A is as same as B. In $D_i$, the second '$L$' presents the prime number. The number of 'R' before this '$L$' is even. So A is an 'even sequence'. In $D_{i+1}$, the second '$L$' in $D_i$ is 'erased' and replaced by 'R'. So B is an 'odd sequence'. We get:

$B > A, D_{i+1} > D_i$.

**End**

According to the MSS theory in symbolic dynamics, we also have:
**Theorem 5:**

$$D_1 < D_2 < D_3 < \cdots < D_i < D_{i+1} < \cdots < D_\infty \\ \Leftrightarrow u(D_1) < u(D_2) < u(D_3) < \cdots < u(D_i) < u(D_{i+1}) < \cdots < u(D_\infty).$$

By elementary number theory, we could prove the period of $D_i$ is $T(D_i) = p_1 \times p_2 \times p_3 \times \cdots \times p_i$. So the period of $D_\infty$ is $T(D_\infty) = p_1 \times p_2 \times p_3 \times \cdots \to \infty$. It seems that prime system will transform from period system into a chaos system. The $D_i, i = 1 \to \infty$ could describe this process.

The relation between symbolic dynamics and one dimension unimodal mapping has been well investigated. So we could use a simple unimodal mapping to describe the dynamics of $D_i, i = 1 \to \infty$. Here we select $x_{n+1} = 1 - ux_n^2$. The relation of its parameter u and related symbols is shown in Table 1.

According to the period expression of prime number, $M_2, M_3, \cdots M_{pi}, \cdots$ are all kneading sequence, which means we could find their corresponding actual orbit from mapping $x_{n+1} = 1 - ux_n^2$. We could also prove their '•' composition $D_1 < D_2 < D_3 < \cdots < D_i < \cdots < D_\infty$ are kneading sequence. So prime system's special road to the chaos could be described by a series of orbits of $x_{n+1} = 1 - u(D_i)x_n^2, u(D_1) < u(D_2) < \cdots < u(D_i) < \cdots < u(D_\infty)$.

For $D_2 = M_2 \bullet M_3 = RLRRRL$, we could get $u(D_2) \approx 1.476$. It's very difficult to get the expression of $D_\infty$ by '•' composition rule. But according to out 'sieve method', all the nature number will be 'erased' except 1. We could get the expression:

$$\lim_{n\to\infty} D_n = RLRRR\cdots = RLR^\infty$$

The mapping $x_{n+1} = 1 - u(D_\infty)x_n^2$ could represent the dynamical character of whole prime system. Here $u(D_\infty) = u(RLR^\infty) \approx 1.5437$.

The Lyapunov exponent of $x_{n+1} = 1 - u(D_\infty)x_n^2$ is about 0.3406. So we could say the primes system is a chaos system. According to the kneading theory, the topological entropy of $D_\infty = RLR^\infty$ is $\ln(2)/2$.

Here we are still not very clear about the approaching method of $D_i, i = 1 \to \infty$ and the chaos orbits near $D_\infty$. So it's just a heuristic proof for Theme1. We will discuss this topic in the further.

## 5 Summary

This paper builds a novel dynamic model for primes system by symbolic dynamics. Then we connect pattern of primes with chaos orbits of Logistic mapping. This bridge may deduce the new prove avenue for twin primes conjecture.